\documentclass[11pt]{article}

\usepackage{a4wide}

\def\versiondate{31. Januar 2006}
\usepackage{amssymb}
\usepackage{amsmath}
\usepackage{amsthm}

\def\comment#1{}

\def\ni{\noindent}

\def\setm{\setminus}

\def\term#1{{\em #1}\marginpar{\raggedright{\small\it #1}}}
\def\nct#1{{\em #1}} %% no-cite term

\def\ITEMMACRO #1 ??? #2 ???{\par\vskip4pt\noindent%
\hangindent=#2em\setbox0\hbox{#1\kern4pt}%
\ifdim\wd0<\hangindent\setbox0\hbox to\hangindent{\hss#1\kern7pt}\fi%
\box0\ignorespaces}
 
\def\Item(#1){\ITEMMACRO {\rm (#1)} ??? 1.8 ???}
\def\IItem#1{\ITEMMACRO #1 ??? 4 ???}
\def\FreeItem#1{\ITEMMACRO #1 ??? 1.8 ???}

\let\Bitem=\bItem
\def\BrackItem[#1]{\ITEMMACRO [#1] ??? 1.8 ???}

\usepackage{graphicx}
\usepackage{color}
\graphicspath{{Figures/}}
\input{psfig.sty}

\newtheorem{theorem}{Theorem}
\newtheorem{lemma}{Lemma}
\newtheorem{defin}{Definition}
\newtheorem{proposition}{Proposition}
\newtheorem{prop}{Proposition}

\newtheorem{problem}{Problem}
\newtheorem{corollary}{Corollary}

\def\Claim#1.{\medbreak\ni{\bf Claim~#1.}}
\def\Case#1.{\medbreak\ni{\bf Case~#1.}}
\def\SubCase#1.{\medbreak\ni{\bf Subcase~#1.}}

\def\Proof{\ni{\sl Proof.}\ }
\def\qed{\hfill\fbox{\hbox{}}\medskip}

\def\below#1{[\kern1pt#1-1\kern-1.5pt\updownarrow\kern-1.5pt#1\kern1pt]}
\def\between#1#2{[\kern1pt#1\kern-1.5pt\updownarrow\kern-1.5pt#2\kern1pt]}
\def\bbelow#1{[\kern1pt#1-1\kern-1.5pt\updownarrow\kern-1pt\updownarrow\kern-1.5pt#1\kern1pt]}
\def\bbetween#1#2{[\kern1pt#1\kern-1.5pt\updownarrow\kern-1pt\updownarrow\kern-1.5pt#2\kern1pt]}

\def\rev#1{#1^{{\rm rev}}}

\def\with{\kern1.5pt :\kern1.5pt}

\def\join{{\vee}}

\def\<{\langle}
\def\>{\rangle}

\def\RR{\hbox{\sf I\kern-1ptR}}
\def\NN{\hbox{\sf I\kern-1ptN}}
\def\ZZ{\hbox{\sf Z\kern-4ptZ}}

\begin{document}
% *************************************************************
%                   TITLE PAGE
% *************************************************************
%

\begingroup\raggedright
\baselineskip22pt
\LARGE\bf

Orthogonal Surfaces

\endgroup
\vskip20pt

\begingroup
\fontsize{12}{14}%
{\sc Stefan Felsner {\rm and} Sarah Kappes}
\vskip6pt

{\small\it Technische Universit\"at Berlin,  Institut f\"ur Mathematik, MA 6-1,}
\hfill\vskip-1pt
{\small\it Stra\ss{}e des 17. Juni 136, 10623 Berlin, Germany}\hfill\vskip1pt
{\small {\it E-mail:} $\{$felsner,kappes$\}$@math.tu-berlin.de}
\endgroup

\footnotetext{\versiondate}

% *************************************************************

\vskip20pt
\begingroup\fontsize{10}{12}\rm
\noindent{\bf Abstract.}
Orthogonal surfaces are nice mathematical objects which have
interesting connections to various fields, e.g.,
integer programming, monomial ideals and  order dimension.
While orthogonal surfaces in one or two dimensions are rather trivial
already the three dimensional case has a rich structure
with connections to Schnyder woods, planar graphs and 3-polytopes.

Our objective is to detect more of the structure of orthogonal
surfaces in four and higher dimensions. In particular we
are driven by the question which non-generic orthogonal
surfaces have a polytopal structure. 

We study characteristic points and the cp-orders of orthogonal
surfaces, i.e., the dominance orders on the characteristic points. In
the generic case these orders are (almost) face lattices of polytopes.
Examples show that in general cp-orders can lack key properties of
face lattices. We investigate extra requirements which may help to
have cp-orders which are face lattices.

Finally, we turn the focus and ask for the realizability of polytopes
on orthogonal surfaces. There are criteria which prevent large classes
of simplicial polytopes from being realizable. On the other hand
we identify some families of polytopes which can be realized
on orthogonal surfaces.

\vskip10pt\ni
\endgroup

% *************************************************************
\section{Introduction}\label{sec:intro}
% *************************************************************

Subsection~\ref{subs:p+m} is a short survey of previous work and
important problems in the field of orthogonal surfaces. 
Subsection~\ref{subs:n+d} is a collection of basic
definitions and notation. 

Section~\ref{sec:3dim} is a review of
3-dimensional surfaces. We briefly look at the generic case and then
move on to non-generic surfaces. In this still well visualizable case
the distinction between generated and characteristic points
becomes obvious and degeneracies can break the otherwise nice properties.
Rigidity is the extra condition which helps.

Section~\ref{sec:odim} relates orthogonal surfaces and order theory.
We discuss Schnyder's characterization of planar graphs and the
Brightwell-Trotter Theorem in their relation with orthogonal surfaces
and explain how dimension theory can help to prove that certain
polytopes are not representable on orthogonal surfaces.

While the first three section mainly collect and review what
was already known the final two sections contain new material.

With Section~\ref{sec:hdim} we move on to higher dimensions.
Issues of degeneracy and the relation between generated and
characteristic points are analyzed with care. Characteristic
points are also compared to the algebraically interesting syzygy-points
of a surface. The concept of rigidity of an orthogonal surface
is generalized to higher dimensions. Two concrete examples 
show that even in the rigid case cp-orders of 4-dimensional
orthogonal surfaces may lack simple properties required for
face-lattices of polytopes.

Section~\ref{sec:real} deals with realizability of polytopes on
orthogonal surfaces. We present a new realizability criterion for
simplicial polytopes.  Exhaustive computations show that this
criterion works for 2344 out of the 2957 simplicial balls on 9
vertices which are obtained by deleting a facet of a non-realizable
polytope.  In the final subsection we identify some families of
realizable polytopes.

% *************************************************************
\subsection{Previous work, motivation}
\label{subs:p+m}

Orthogonal surfaces have been studied by Scarf~\cite{s-cee-73} in the
context of test sets for integer programs.  Initiated by work of
Bayer, Peeva and Sturmfels~\cite{bps-mr-98} they later became of
interest in commutative algebra. A recent monograph of Miller and
Sturmfels~\cite{ms-cca-04} presents the state of the art in this
area. Miller~\cite{m-pgmrtmi-02} was the first to observe the
connections between orthogonal surfaces, Schnyder woods and the
Brigthwell-Trotter Theorem about the order dimension of polytopes.  We
will outline these connections in Sections~\ref{sec:3dim} and
\ref{sec:odim} where we also review other applications of order
theoretic results to orthogonal surfaces.

Before stating the Theorem of Scarf which can be regarded the most fundamental
result in the field we briefly set the stage with the most important
terms. 

Our starting point is a (finite) antichain $V$ in the dominance order
on $\RR^d$. The \term{orthogonal surface} $S_V$ generated by $V$ is
the topological boundary of the filter  
$\< V \> = \big\{x\in \RR^d : \text{ there is a } v\in V \text{ with } v_i \leq x_i
\text{ for all } i\big\}$.

An orthogonal surface $S_V$ in $\RR^d$ is \term{suspended} 
if $V$ contains $d$ extremal vertices. An orthogonal surface $S_V$ is
\term{generic} if no two points in $V$ have a coordinate in common.

The \term{Scarf complex} $\Delta_V$ of a generic orthogonal surface $S_V$
generated by $V$ consists of all the subsets $U$ of $V$ with the
property that $\bigvee_{u\in U} u \in S_V$. It is a good exercise to show that
$\Delta_V$ is a simplicial complex.

\begin{theorem}[Scarf '73]\label{thm:scarf}
The Scarf complex $\Delta_V$ of a generic suspended orthogonal surface $S_V$
in $\RR^d$ is isomorphic to the face complex of a simplicial
$d$-polytope with one facet removed.
\end{theorem}

A proof of the theorem is given in~\cite{bps-mr-98}. 
Figure~\ref{fig:octahedron} shows an example. 
The dimension 3 case of Scarf's theorem was independently discovered by
Schnyder~\cite{s-pgpd-89}. 

An interesting problem inspired by Scarf's theorem is the realization
question, asking for a characterization of those simplicial
$d$-polytopes which have a corresponding orthogonal surface. We come
back to this question in Sections~\ref{sec:odim} and~\ref{sec:real}.

%%%%%%%%%%%%%%%%%%%%%%%%%%%%%%%%%%%%%%%%%%%%%%%%%%%
%%
\PsFigCap{60}{octahedron}%
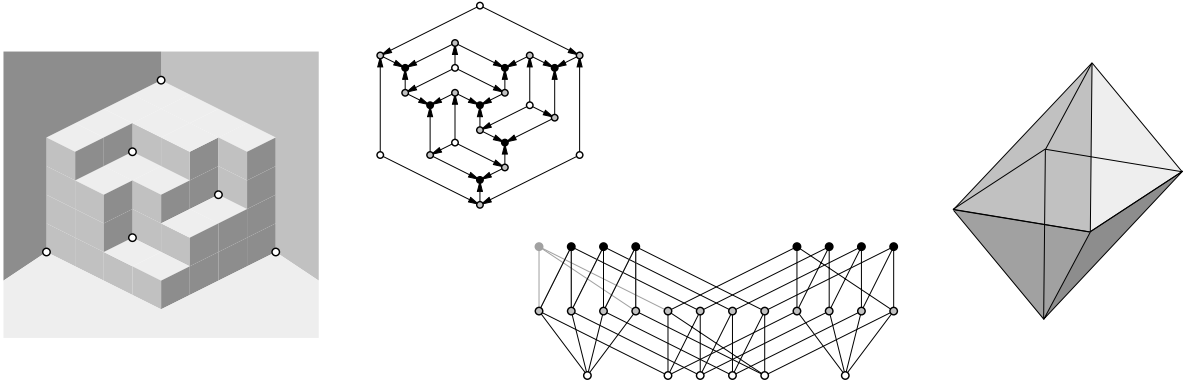
{An orthogonal surface, two diagrams of its complex
and the corresponding polytope.}
%%
%%%%%%%%%%%%%%%%%%%%%%%%%%%%%%%%%%%%%%%%%%%%%%%%%%%

The subject becomes much more complicated if we consider non-generic
surfaces. In this case, it is not even clear how to define an
appropriate complex on the vertex set $V$.
To overcome this difficulty, we introduce an alternative
interpretation of the Scarf-complex. We observe
that every element $U \in \Delta_V$ corresponds to a \term{characteristic
  point} $p_U = \bigvee(U) \in S_V$. A more general
definition of characteristic points is given in Section \ref{sec:3dim}. For
now, it is sufficient to think of them as the corners of the staircase.

The \term{cp-order} of an orthogonal surface is the set of
characteristic points equipped with the dominance order together with
artificial $0$ and $1$ elements.
With this terminology, we can rephrase Scarf's theorem as follows:

\begin{theorem}[Scarf '73]\label{thm:scarf2}
The cp-order of a generic suspended orthogonal surface is isomorphic
to the face-lattice of some simplicial $d$-polytope with one facet
removed.
\end{theorem}

If the cp-order is a lattice, we will call it a
\term{cp-lattice}. Scarf's Theorem implies that this is always the
case if $S_V$ is generic. 

One of our main goals is to determine conditions that
are less restrictive than genericity but still guarantee that the
cp-order has strong properties. In general, cp-orders of non-generic
surfaces are no lattices, not graded and do not satisfy the
diamond-property. We will discuss examples in Sections~\ref{sec:3dim}
and ~\ref{sec:hdim}.
To deal with these problems for 3-dimensional orthogonal surfaces,
Miller introduced the notion of
\term{rigidity}, \cite{m-pgmrtmi-02}. We will define and discuss 
this property in Section~\ref{sec:3dim} and its generalization to
higher dimensions in Section~\ref{sec:hdim}.
 
The following theorem comprises a generalization of Scarf's
theorem and the solution for the realization problem for the 3-dimensional
case:

\begin{theorem}\label{thm:main3d}
The cp-orders of rigid suspended orthogonal surfaces  in $\RR^3$
correspond to the face-lattices of $3$-polytopes with one facet
removed.
\end{theorem}

In particular, rigidity implies that the cp-order is graded and a lattice.
This result can be regarded as a strengthening of the
Brightwell-Trotter Theorem~\cite{bt-odcp-93} about the order dimension
of face lattices of 3-polytopes (Theorem~\ref{thm:B-T}). Proofs can
be found in~\cite{f-gepg-03} and~\cite{fz-os3d-06}. These proofs
actually show more, namely a bijection with Schnyder woods. We review
some aspects of the theory in Section~\ref{sec:3dim}.

For dimensions $d>3$ and non-generic surfaces, it is already
challenging to come up with appropriate combinatorial definitions for
characteristic points and properties of the cp-order. We present some
results in Section \ref{sec:hdim}. The dream which originated this
research was to obtain some high-dimensional generalization of
Theorem~\ref{thm:main3d}. The dream did not become true but we have
shaped some basic blocks of theory which should have future.

%%%%%%%%%%%%%%%%%
\subsection{Basic notation and definitions}
\label{subs:n+d}

We consider $\RR^d$ equipped with the \term{dominance order}, this is the
partial order on the points
defined by the product of the orders of components, i.e. for $v, w \in \RR^d$:
$$
(v_1, v_2, \dots, v_d) \leq (w_1, \dots, w_d)
\quad\Longleftrightarrow\quad
v_i \leq w_i \text{ for all } i \in \{1, \dots,d\}.
$$

We say that $v$ \term{strictly dominates} $w$ if $v_i \gneq w_i$ for
all $i = 1, \dots, d$ and denote this relation by $v \rhd w$. 

A point $v$ \term{almost strictly dominates}
another point $w$, if $v_i = w_i$ for exactly one coordinate $i$ and
$v_j \gneq w_j$ for all $j \neq i$, we denote this with
$v \rhd_i w$.

The \term{join} $v \vee w$ of points $v$ and $w$ is defined as the
componentwise maximum of $v$ and $w$ and the \term{meet} $v \wedge w$
as the componentwise minimum of $v$ and $w$.

The \term{cone} $C(v)$ of $v\in\RR^d$ is  the set of all points greater
than $v$ in the dominance order, formally
$C(v) = \{ x \in \RR^d \text{ }|\text{ } x \geq v\}$.

An \term{antichain} $V \subset \RR^d$ is a set
of pairwise incomparable points. This means for any $v, w \in  V$,
there are two coordinates $i, j \in \{1, \dots, d\}$ such that $v_i <
w_i$ and $v_j > w_j$. Equivalently, no point of $V$ is contained in
the cone of any other.

The \term{filter} $\< V \>$ generated by $V$ is the union of
all cones $C(v)$ for $v \in V$.

The \term{orthogonal surface} $S_V$ generated by $V$ is
the boundary of $\langle V \rangle$. The generating set $V$ is an antichain
exactly if all elements of $V$ appear as minima on $S_V$.
We will generally assume that this is the case. 

A point $p \in \RR^d$ belongs to $S_V$ if and only if there is
a vertex $v \in V$ such that $v \leq p$ and there is \emph{no} $w \in V$ such
that $p \rhd w$. In other words, $S_V$
consists of points that share some coordinate with every vertex $v\in V$ they
dominate.

With a point $p \in S_V$, we associate a \term{down-set} $D_p =\{v \in
V : v \leq p \}$. For a point $p$ and $v \in D_p$ define $T_p(v) = \{i
\in \{1, \dots, d\}: p_i = v_i\}$, this is the set of \term{tight
coordinates} of $p$ with respect to $v$.

An orthogonal surface $S_V$ in $\RR^d$ is \term{suspended} if $V$
contains a suspension vertex for each $i$, i.e., a vertex with
coordinates $(0,\ldots,0,M_i,0,\ldots,0)$, for each $i$ and the
coordinates of each non-suspension vertex $v\in V$ satisfy $0\leq v_i
< M_i$.

An orthogonal surface $S_V$ is \term{generic} if no two
points in $V$ have the same $i$th coordinate, for any $i$.
If $S_V$ is suspended then the condition has to be relaxed for the
suspension vertices which obviously share coordinates of value zero.

The \term{Scarf complex} $\Delta_V$ of a generic orthogonal surface $S_V$
generated by $V$ consists of all the subsets $U$ of $V$ with the
property that $\bigvee_{u\in U} u \in S_V$. 

Figure~\ref{fig:examp} shows orthogonal surfaces in two and three
dimensions.  The picture on the right is obtained as orthogonal
projection onto the plane \mbox{$x_1+x_2+x_3 = 0$}.

%%%%%%%%%%%%%%%%%%%%%%%%%%%%%%%%%%%%%%%%%%%%%%%%%%%
%%
\PsFigCap{40}{examp}%%
{Orthogonal surfaces in two and three dimensions.}
%%
%%%%%%%%%%%%%%%%%%%%%%%%%%%%%%%%%%%%%%%%%%%%%%%%%%%

% *************************************************************
\section{The 3-dimensional case}\label{sec:3dim}
% *************************************************************

In this section we discuss 3-dimensional orthogonal surfaces. 
Before turning to the general case it seems appropriate
to review the main correspondence in the generic case.
We start with a correspondence between characteristic points 
of the surface and elements of the Scarf-complex:

\Bitem Rank $0$ elements of the complex (vertices) correspond to the
minima of the surface, i.e, to elements of $V$.
\Bitem Rank $1$ elements of the complex (edges) correspond to those
elements of the surface which can be written as join $u\join v$
for a pair $u,v$ of vertices.
\Bitem Rank $2$ elements of the complex (faces) correspond to the
maxima of the surface, alternatively, these elements are joins
of triples of vertices.

% *************************************************************
\subsection{3-Dimensional and generic}
% *************************************************************

Let $S_V$ be a generic suspended orthogonal surface in $\RR^3$, i.e., no two
non-suspension points in $V$ have a coordinate in common.
We identify the coordinates 1,2,3 with the colors red, green and blue, in
this order. In addition, we assume a cyclic structure on the coordinates
such that $i+1$ and $i-1$ is always defined.

It is valuable to have a notation for some special features of the
surface. For a vertex $v\in V$ and a color $i$ define the 
\term{flat}\footnote{The definition given here is only
valid in the generic case!}
$F_i(v)$ as the set of points on $S_V$ which dominate $v$ and share
coordinate $i$ with $v$. The intersection $F_{i-1}(v) \cap F_{i+1}(v)$
of two flats of $v$ is the \term{orthogonal arc} of $v$ in color $i$.

Draw every rank 1 element $\{u,v\}$ of the complex $\Delta_V$ as combination of two
straight line segments, one connecting $u$ to $u\join v$, the
other connecting $v$ to $u\join v$. This yields 
a drawing of a graph on $S_V$. Before discussing
properties of the graph we impose additional structure on these edges.

For two vertices $u$ and $v$ (at most one of them a suspension vertex)
genericity implies that the join $u\join v$ has one coordinate
from one of the them and two of the coordinates from the other
vertex. In particular this is true if $u\join v \in S_V$, i.e, if
$u,v$ is an edge in $\Delta_V$.  If $u\join v \in S_V$ and $u\join v$
has two coordinates from $v$ then we orient the edge as $v\to u$ and
color it with color of the coordinate which comes from $u$. In
Figure \ref{fig:edge-ex} $u\join v = (v_1,v_2,u_3)$ and the edge is oriented $v\to
u$ and colored with color 3. The drawn edge consists of the
orthogonal arc of $v$ in color 3 which leads from $v$ to $u\join v$ 
and a segment between $u\join v$ and $u$ which traverses the flat
$F_3(u)$.

%%%%%%%%%%%%%%%%%%%%%%%%%%%%%%%%%%%%%%%%%%%%%%%%%%%
%%
\PsFigCap{80}{edge-ex}% 
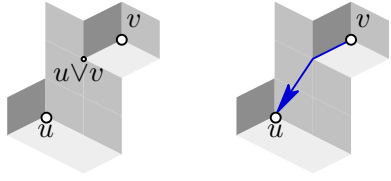
{Drawing and orienting edges on $S_V$.}
%%
%%%%%%%%%%%%%%%%%%%%%%%%%%%%%%%%%%%%%%%%%%%%%%%%%%%

From the geometry of flats (in the generic case a flat only contains
a single vertex $v\in V$) and the way edges are drawn we can conclude
the following:

\Bitem There are no crossing edges, i.e, the graph is planar.

\Bitem Every maximum of the surface dominates exactly three vertices,
i.e., the graph is a triangulation.

\Bitem The orientation and coloring of the edges has the following
properties:
\IItem{\qquad[\;Rule of Vertices\;] } Every non-suspension vertex $v$ has one
outgoing edge in each color. The out-edges $e_1, e_2, e_3$ with colors
1, 2, 3 leave $v$ in clockwise order. Each edge entering $v$ with
color $i$ enters in the clockwise section from $e_{i+1}$ to $e_{i-1}$.
Suspension vertices only have incoming edges of one color.

\FreeItem{} The `rule of vertices'  defines a \term{Schnyder wood}
of a planar triangulation.
\medskip

This explains how to obtain a Schnyder wood on a planar
triangulation from a generic suspended orthogonal surface in
$\RR^3$. For the converse consider a triangulated planar
graph. Selecting an outer triangle yields an essentially unique plane
embedding.  Specify a Schnyder wood of the plane triangulation -- it was shown 
by Schnyder~\cite{s-pgpd-89} that these structures exist,
actually, a triangulation can have many different Schnyder woods, see
\cite{f-lspg-04}.

The set of all edges of color $i$ forms a directed tree spanning all
interior vertices of the triangulation, this tree is rooted at one of
the three outer vertices which will be called the \term{suspension
vertex} of color $i$. The three trees define three colored paths
$P_1(v)$, $P_2(v)$ and $P_2(v)$ from an interior vertex $v$ to the
three outer vertices. From the `rule of vertices' it can be deduced
that these paths are interiorly disjoint. Hence, they partition the
interior of the outer triangle of the graph into three regions
$R_1(v)$, $R_2(v)$ and $R_3(v)$. The \term{region vector} of a vertex
$v$ is the vector $(v_1,v_2,v_3)$ defined by 

$$ v_i = \hbox{The number of faces contained in region $R_i(v)$}.  $$

The set of region vectors of vertices of the graph yields a finite
antichain $V\subset \RR^3$ such that the orthogonal surface $S_V$ has
a complex $\Delta_V$ which is isomorphic to the original plane triangulation.
Moreover, the orientation and coloring of edges on the surface $S_V$
induced the Schnyder wood used for the construction of the surface.

Some of the details of the proof can be found in the original papers
of Schnyder~\cite{s-pgpd-89,s-epgg-90}, the notion of an orthogonal surface,
however, was not known to Schnyder. Proofs given in the
publications \cite{m-pgmrtmi-02,f-gepg-03} and in the book 
\cite{f-gga-04} extend these ideas to the more
general case of 3-connected planar graphs.

By Steinitz's Theorem planar triangulations are essentially the same
as simplicial 3-polytopes. Disregarding a facet of a 3-polytope corresponds to a
choice of the outer face for the corresponding planar graph.
Therefore, the following proposition is a colored strengthening of special cases 
of Theorem~\ref{thm:scarf} and Theorem~\ref{thm:main3d}:

\begin{proposition}
The edge colored complex of a generic suspended orthogonal surface
$S_V$ in $\RR^3$ is the Schnyder wood of a plane
triangulation. Moreover, every Schnyder wood of a plane triangulation
has a corresponding orthogonal surface.
\end{proposition}

In the above sketch we have been using Schnyder
woods. In~\cite{s-epgg-90} Schnyder introduced \term{angle labelings}
of plane triangulations and proved that they are in bijection with
Schnyder woods. The two properties of angle labelings are

\BrackItem[\;Rule of Vertices\;] The labels of the angles at each vertex form,
in clockwise order, a non-empty interval of 1's, 2's and 3's.
\BrackItem[\;Rule of Faces\;] The labels in each face are 1, 2, 3 in clockwise order.
\medskip

\ni
From an orthogonal surface supporting a Schnyder wood the corresponding
angle labeling is directly visible: The angle between consecutive edges $e$ and 
$e'$ at vertex $v$ is colored~$i$ if both edges leave $v$ on the flat $F_i(v)$.
Figure~\ref{fig:dict-1} shows a graph on a surface with the induced
edge and angle colorings.

%%%%%%%%%%%%%%%%%%%%%%%%%%%%%%%%%%%%%%%%%%%%%%%%%%%
%%
{\def\SetFigFont#1#2#3#4#5{\tiny} %% Fuer das Bild!!
\PsFigCap{65}{dict-1}%
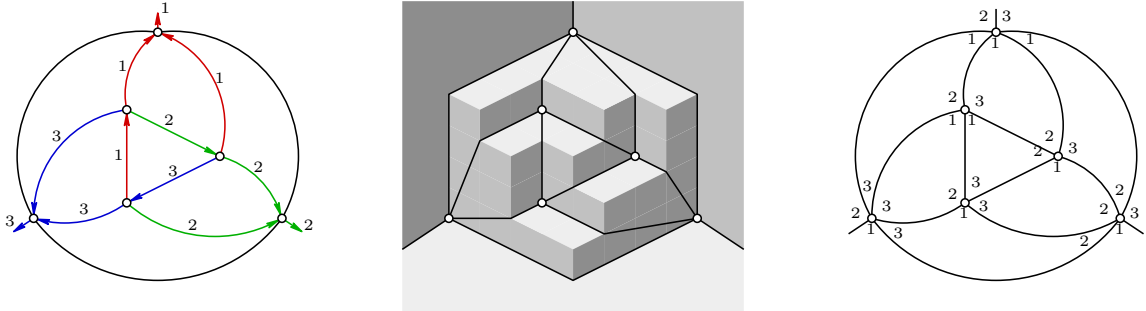
{The graph of an orthogonal surface and corresponding Schnyder's colorings.}
}
%%
%%%%%%%%%%%%%%%%%%%%%%%%%%%%%%%%%%%%%%%%%%%%%%%%%%%

% *************************************************************
\subsection{3-Dimensional and non-generic}
% *************************************************************

\ni Given a non-generic antichain $V$ in $\RR^3$ it would be nice to have
a complex $\Delta_V$ such that the elements of the complex are
in bijection with the characteristic points of the surface $S_V$,
just as in the generic case.
Attempts to define such a $\Delta_V$ face some problems.

First of all, we have to rework and generalize our notion of a
flat. Instead of attaching a flat strictly to one minimum, we now
think of flats as connected $(d-1)$-dimensional components of the
intersection of $S_V$ with some hyperplane. In the non-generic case
such a component/flat can contain several minima, all sharing the coordinate,
which defines the flat, see Figure~\ref{fig:flats}.

%%%%%%%%%%%%%%%%%%%%%%%%%%%%%%%%%%%%%%%%%%%%%%%%%%%
%%
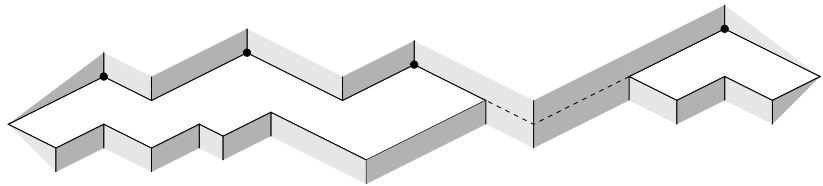
\PsFigCap{50}{flats}{Two flats with the same defining coordinate, one with three
  minima and one with a single minimum.}
%%
%%%%%%%%%%%%%%%%%%%%%%%%%%%%%%%%%%%%%%%%%%%%%%%%%%%

For every $v \in V$ and every coordinate $i$, the
\term{almost strict upset} $U_i(v) = \{p \in S_V: p \rhd_i v\}$
belongs to the same $i$-flat as $v$. 

If $U_i(v) \cap U_i(w) \neq \emptyset$, then $v$ and $w$ belong to
the same $i$-flat. More general, we define a relation $\thicksim_i$ on
$V$ by $v \thicksim_i w \Leftrightarrow U_i(v) \cap U_i(w) \neq
\emptyset$. The transitive closure $\thicksim_i^c$ of $\thicksim_i$ is
an equivalence relation. The equivalence classes are exactly those
sets of minima sharing a common $i$-flat.

\begin{defin}\label{def:flat} 
Let $v \in V$. The \term{$i$-flat} $F_i(v)$ is
the topological closure of the set $$ \bigcup \limits_{w \thicksim^c_i
v} U_i(w)$$ \end{defin}

The equivalence class of minima on an $i$-flat $F_i$ is $V_{F_i} = F_i
\cap V$.  Furthermore, we define the \term{upper part of the flat}
$F_i$ as $F_i^u = \bigcup \limits_{v \in V_{F_i}} U_i(v) = \{ p \in
S_V: p \rhd_i v \text{ for some } v \in V_{F_i} \}$

% *************************************************************
\subsubsection{Degeneracies}

There may be characteristic points which can be obtained as join of
distinct pairs of vertices, e.g., $p = u\join v = v\join w = w\join u$
for distinct vertices $u,v,w$. Figure~\ref{fig:degenerate} shows an example. We
want to have the property that every orthogonal arc is part of an
edge, i.e., connects a vertex with a characteristic point of
rank~1. Therfore, we usually assume that surfaces in $\RR^3$
have no such substructure, if we want to emphasize this property we say
the surface is \term{non-degenerate}.

%%%%%%%%%%%%%%%%%%%%%%%%%%%%%%%%%%%%%%%%%%%%%%%%%%%
%%
\PsFigCap{35}{degenerate}%
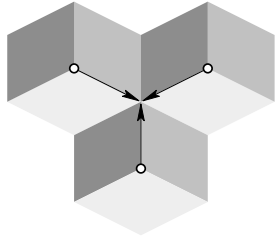
{A degenerate situation on an orthogonal surface.}
%%
%%%%%%%%%%%%%%%%%%%%%%%%%%%%%%%%%%%%%%%%%%%%%%%%%%%

% *************************************************************
\subsubsection{Generated versus characteristic points}

In the generic case every joint $u\join v$ on the surface $S_V$ is
a characteristic point and corresponds to a rank~1 element of
$\Delta_V$. This is not true in general. An example is shown in
Figure~\ref{fig:charAndNonChar} where  characteristic points are black
but there are additional (white) generated points.

%%%%%%%%%%%%%%%%%%%%%%%%%%%%%%%%%%%%%%%%%%%%%%%%%%%
%%
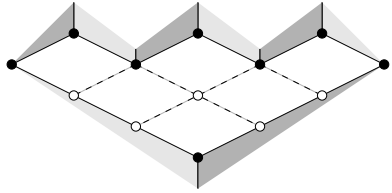
\PsFigCap{65}{charAndNonChar}{Characteristic and non-characteristic joins}
%%
%%%%%%%%%%%%%%%%%%%%%%%%%%%%%%%%%%%%%%%%%%%%%%%%%%%

This shows the need of a new definition for characteristic points.
In dimension three we could stick to the definition that
characteristic points of rank~1 are endpoints of orthogonal arcs
while all other characteristic points are minima or maxima.
More satisfactory and more appropriate for generalizations to higher dimensions
is the following:

\begin{defin}\label{def:c-point} 
A \term{characteristic point} is a point which is incident to flats of all 
colors.
\end{defin}

Clearly, every minimum is a characteristic point. From the definition
it is immediate that every characteristic point is a generated point, i.e. 
can be expressed as the join of some minima.
Figure~\ref{fig:typen} shows the possible
types of points on a surface. Characteristic points
are those of types {\bf a}, {\bf d} and {\bf e}. Type {\bf e} is the
forbidden degenerate substructure.

%%%%%%%%%%%%%%%%%%%%%%%%%%%%%%%%%%%%%%%%%%%%%%%%%%%
%%
%\PsTexFig{flat}%
%{The flats of the intersection of $S_V$ with $x_1=1$.}
%%\centerline{FIG flat}
%%
%%%%%%%%%%%%%%%%%%%%%%%%%%%%%%%%%%%%%%%%%%%%%%%%%%%

%%%%%%%%%%%%%%%%%%%%%%%%%%%%%%%%%%%%%%%%%%%%%%%%%%%%%%%%%%%%%%%
%%
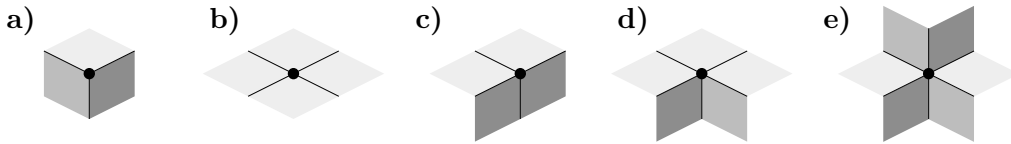
\PsFigCap{95}{typen}{A classification of point types on orthogonal surfaces in $\RR^3$.}
%%
%%%%%%%%%%%%%%%%%%%%%%%%%%%%%%%%%%%%%%%%%%%%%%%%%%%%%%%%%%%%%%%

% *************************************************************
\subsubsection{Rigidity}

In the generic case the dominance order on characteristic points and the
inclusion order of the sets of the complex $\Delta_V$ coincide.
This is no longer true in the general case. Even characteristic point of rank~1
can dominate many vertices, see Figure~\ref{fig:non-rigid}.

%%%%%%%%%%%%%%%%%%%%%%%%%%%%%%%%%%%%%%%%%%%%%%%%%%%%%%%%%%%%%%%
%%
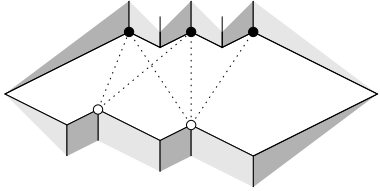
\PsFigCap{65}{non-rigid}{Characteristic points of rank~1 with 
 two and three dominated vertices.}
%%
%%%%%%%%%%%%%%%%%%%%%%%%%%%%%%%%%%%%%%%%%%%%%%%%%%%%%%%%%%%%%%%

In this case the graph defined by the surface is not unique.
Uncoordinated choices for edges can even lead to crossing edges.
Miller~\cite{m-pgmrtmi-02} calls a surface \term{rigid} if characteristic
points $u\join v$ only dominate $u$ and $v$ in $V$, i.e., there is no
$w\in V\setminus\{u,v\}$ with $w \leq u\join v$.

Note that rigidity of a surface $S_V$ in $\RR^3$ implies that $S_V$
is non-degenerate and that -as in the generic case- $S_V$ defines a
unique graph on the vertex set $V$: $(u,v)$ is an edge if and only if
$u \join v$ is a characteristic point. 
Such an edge can be drawn as the combination of two straight line
segments $u - u\join v$ and  $v - u \join v$. 
At least one of them is an orthogonal arc of the surface and all
orthogonal arcs emanating from vertices $v\in V$ are used.
Again there is an obvious definition of orientations and
colorings of edges.

Let $G$ be a plane graph with 
\term{suspension vertices} $a_1,a_2,a_3$ on the outer face.
Add a half-edge to each of the three suspension vertices.
A \term{Schnyder wood} for $G$ is an orientation and coloring
of the edges with colors $1,2,3$
such that:

\Item(W1)
Every edge $e$ is oriented by one or two opposite directions.
The directions of edges are labeled such that 
if $e$ is bioriented, then the two directions have
distinct labels.

\Item(W2)
The half-edge at $a_i$ is directed outwards and labeled $i$.

\Item(W3)
Every vertex $v$ has outdegree one in each label.
The edges $e_1,e_2,e_3$ leaving $v$ in labels 1,2,3 occur in clockwise
order. Each edge entering $v$ in label $i$ enters $v$ in the clockwise
sector from $e_{i+1}$ to $e_{i-1}$. [\; Rule of vertices\;]

\Item(W4)
There is no interior face whose boundary is a directed cycle
in one label. 
\medskip

\ni
The orientation and coloring of edges induced by a suspended rigid orthogonal surface
is a Schnyder wood for the induced plane graph. 

It can be shown that a plane graph $G$ with suspension vertices
$a_1,a_2,a_3$ has a Schnyder wood exactly if the graph $G_\infty$
which is obtained from $G$ by adding a new vertex adjacent to the
three suspension vertices is 3-connected.

Let $G$ be a plane graph with a Schnyder wood. The edges of color $i$
in the Schnyder wood induce a spanning tree rooted at $a_i$. These
trees define paths and these paths, in turn, define three regions for
every vertex. Therfore we can again consider the set $V$ of region
vectors of the vertices. This set $V$ is an antichain in $\RR^3$ and
the surface $S_V$ generated by $V$ supports the graph $G$ and the
Schnyder wood which was used to define the regions.
However, the surface $S_V$ obtained by this construction needs not to be rigid.

Let $S_V$ be constructed from the region vectors of a graph $G$ with a
Schnyder wood. The vertices, edges and bounded faces of $G$
can be associated to the characteristic points of $S_V$.
Actually, it is even possible to associate with a Schnyder wood on a
plane graph a {\em rigid} orthogonal surface which supports the 
Schnyder wood. Hence, there is an orthogonal surface which uniquely
supports the given Schnyder wood. This has been conjectured by
Miller~\cite{m-pgmrtmi-02} and was proven in~\cite{f-gepg-03}
and~\cite{fz-os3d-06}, we come back to this in the next section.

With Steinitz's correspondence between 3-connected planar graphs and
3-polytopes we obtain Theorem~\ref{thm:main3d-2}. This theorem is a
more precise restatement of Theorem~\ref{thm:main3d}.

\begin{theorem}\label{thm:main3d-2}
The cp-orders of rigid suspended orthogonal surfaces 
in $\RR^3$ coincide with the face-lattices of 
3-polytopes with one facet removed or with one vertex of degree three
and the incident faces removed.
\end{theorem}

% *************************************************************
\subsubsection{Duality}

Let $V$ be the generating antichain for a rigid orthogonal surface
$S_V$ in $\RR^3$ and let $W$ be the set of maxima of $S_V$.  Consider
the reflection at $\bf 0$ and let $\bar{W}$ be the image of $W$ under
this map. The orthogonal surface $S_{\bar{W}}$ turns out to be
-almost\footnote{The difference between $S_V$ and $S_{\bar{W}}$ is in
the unbounded flats.}- the same surface as $S_V$ with the reversed
direction of the dominance order.  The surface $S_{\bar{W}}$ is again
rigid and supports a unique Schnyder wood. Well not quite, the surface
$S_{\bar{W}}$ can have more than three unbounded orthogonal arcs. This
can be repaired by adding three suspension vertices to the dual which
bundle the unbounded orthogonal arcs. Figure~\ref{fig:dual-os} shows
an example.  A more detailed account to the duality of Schnyder woods
can be found in~\cite{f-lspg-04}.

%%%%%%%%%%%%%%%%%%%%%%%%%%%%%%%%%%%%%%%%%%%%%%%%%%%%%%%%%%%%%%%
%%
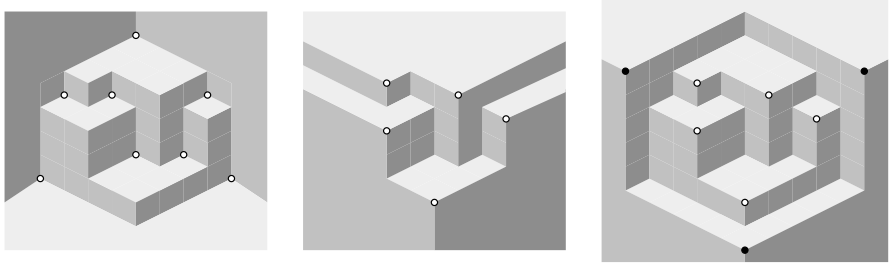
\PsFigCap{50}{dual-os}{A rigid surface, the dual surface
                       and the suspended dual surface.}
%%
%%%%%%%%%%%%%%%%%%%%%%%%%%%%%%%%%%%%%%%%%%%%%%%%%%%%%%%%%%%%%%%

One interesting aspect of the duality is that superimposing a Schnyder
wood and its dual Schnyder wood induces a decomposition of the surface
into quadrangular patches. Each of these patches is completely
contained in a flat, i.e., we can associate a color with each patch.
This yields a joint angle coloring of the underlying planar graph and
the dual.

% *************************************************************
\section{Orthogonal surfaces and order dimension}\label{sec:odim}
% *************************************************************

Every order $P=(X,\leq)$ can be represented as intersection of
linear extensions. That is there are linear orders
$L_1,\ldots,L_k$ with the following properties:
\Bitem If $x \leq y$ in $P$, then $x \leq y$ in each $L_i$,
i.e., the $L_i$ are linear extensions of $P$.
\Bitem If $x||y$ in $P$, then there are indices $i$ and $j$ such that
$x < y$ in $L_i$ and $y < x$ in each $L_j$.
\medskip

\ni
A set of linear extensions representing $P$ in this sense is 
called a \term{realizer} of $P$.
The smallest  number of linear extensions in a realizer of $P$ is the
\term{dimension} $\dim(P)$ of $P$. 

Let $L_1,\ldots,L_k$ be a realizer of $P$. With every $x\in X$
associate a vector $(x_1,\ldots,x_k)\in \RR^k$, where $x_i$ gives
the position (coordinate) of $x$ in $L_i$. This mapping of the
elements of $P$ to points of $\RR^k$ embeds $P$ into the dominance
order of $\RR^k$. Ore defined $\dim(P)$ as the minimum~$k$ such that
$P$ embeds into $\RR^k$ in this way. To prove that the two
definitions are equivalent it remains to show how to obtain a
realizer from an order preserving embedding into $\RR^k$. If the
coordinates of all points in the embedding are pairwise different 
(general position), then the
projections to the coordinate axes form a realizer. Otherwise let
$Y\subset X$ be a set of points sharing a coordinate, e.g., $Y =
\{x\in X : x_1=a \}$.  The order relation among points in $Y$ is
completely determined by the coordinates $2,..,k$. Let $L(Y)$ be any
linear extension of the order induced by $Y$.  Displace coordinate 1
of the elements in $Y$ by tiny amounts such that their projection
confines with $L(Y)$. Repeated perturbations of this type yield an embedding
of $P$ in $\RR^k$ which is in general position.

The following proposition is evident:
\begin{proposition}\label{prop:dim}
Let $X$ be a finite set of points on an orthogonal surface $S_V$ in $\RR^k$
and let $P=(X,\leq)$ be the dominance order on $X$, then $\dim(P) \leq k$.
\end{proposition}

With a graph $G=(V,E)$ associate the incidence order $P_G$ as the order
on the set $V\cup E$ with relations $v \leq e$ iff $v$ is one of the
two vertices of $e$. Schnyder's celebrated characterization of planar graphs
is the following:

\begin{theorem}[Schnyder]
A graph $G$ is planar iff $\dim(P_G) \leq 3$.
\end{theorem}

It was known already to Babai and Duffus~\cite{bd-dagl-81} that $\dim(P_G)\leq 3$
implies that $G$ is planar. Schnyder contributed the other direction.
A proof in our context can follow these steps:
Add edges to $G$ to produce a planar triangulation $G^*$. Using a Schnyder
wood this triangulation can be embedded in a generic orthogonal surface $S_V$
in $\RR^3$. The dominance order on the characteristic points of $S_V$ is
isomorphic to the complex $\Delta_V$ which contains $P_G$ as a suborder.
From that the result follows with Proposition~\ref{prop:dim}.
 
Actually, the above sketch shows that for planar triangulations the
incidence order of vertices, edges and bounded faces has dimension at
most 3. This was known to Schnyder (see~\cite{s-pgpd-89}), it is the
simplicial polytope case of the following generalization of
Schnyder's Theorem.

\begin{theorem}[Brightwell-Trotter] \label{thm:B-T}
Let $P$ be the inclusion order of vertices,
 and faces of a 3-polytope, then $\dim(P) =4$.
The inclusion order of vertices, edges
and all but one of the faces only has dimension 3.
\end{theorem}

The first part is based on a lower bound for the dimension of face lattices
of polytopes. If $P$ is a $d$-polytope and ${\cal F}(P)$ is its face lattice,
then $\dim({\cal F}(P)) \geq d+1$. Since all the critical pairs are between
a maximal and a minimal element of ${\cal F}(P)$ the bound on the dimension  
already holds for the suborder induced by maximal and a minimal elements.

The second part follows from the existence of a rigid embedding of the corresponding
graph on an orthogonal surface in $\RR^3$.

% *************************************************************
\subsubsection{Realizability and order theory}
\label{ssec:rot}

In the terminology developed in the meanwhile we can restate Scarf's
theorem: The dominance order on characteristic points of a generic suspended
surface in $\RR^d$ is isomorphic to the face lattice of a simplicial
$d$-polytope with one facet removed. This result motivates the
following general question:

\begin{problem}[Realizability Problem] 
Which $d$-polytopes can be \term{realized}
on an orthogonal surface in $\RR^d$, i.e., which face lattices of 
$d$-polytopes, with one facet removed are cp-lattices?
\end{problem}

Order theory can provide some criteria for non-realizability.
\medskip

The \term{dimension of the complete graph} $K_n$ is the dimension of 
its incidence order. The asymptotic behavior of this parameter
was first discussed by Spencer~\cite{s-mssso-72}. Trotter improved the lower
bound. Their work implied that the dimension of the complete graph
is closely related to
the number of antichains in the subset lattice. This well studied
problem is known as ``Dedekind's Problem.''  Although no closed form
answer is known, good asymptotic bounds are known, they suffice to show that
$$
\dim(K_n)\sim \log\log n +\bigl(1/2+o(1)\bigr)\log\log\log n.
$$
More recently Ho\c{s}ten and Morris~\cite{hm-odcg-99} could directly relate
$\dim(K_n)$ to a specific class of antichains in the subset lattice.
From this work we know the precise value of $\dim(K_n)$ for all $n\leq 10^{20}$,
for example $\dim(K_{12}) = 4$ and $\dim(K_{13}) = 5$.

For all integers $n$ there exist simplicial 4-polytopes with a
complete graph as skeleton, i.e., the first two levels of their face
lattice is the incidence order of a complete graph $K_n$. These
polytopes are called neighborly (c.f. Ziegler~\cite{z-lp-94}).
From the dimension of complete graphs it follows
that for $n\geq 13$ these 4-polytopes are not realizable on an orthogonal surface in
$\RR^4$.

A more general criterion was developed by Agnarsson, Felsner and
Trotter~\cite{aft-mnegb-99}. They show that the number of edges of a
graph with an incidence order of dimension 4 can be at most $\frac{3}{8}n^2 + o(n^2)$.

With increasing dimension $d$ there is only a rather weak
bound: From $\dim(K_r) > d$ it can be concluded that a graph of
dimension $d$ has at most $\frac{1}{2}(1-\frac{1}{r}) n^2$ edges.
For $d=5$ this gives  a bound of $\frac{81}{164} n^2$ edges.
\medskip

Orthogonal surfaces are completely determined by the position of their vertices.
Therefore, the following notion for the \term{dimension
of a graph}, seems to be more
appropriate in our context.
Let $G=(V,E)$ be a finite simple graph.  A nonempty family
$\mathcal{R}$ of linear orders on the vertex set $V$ of graph $G$ is
called a \term{realizer} of $G$ provided:

\Item($*$) For every edge $e\in E$ and every vertex $x\in V\setm e$,
there is some $L\in\mathcal{R}$ so that $x>y$ in $L$ for
every $y\in e$.

\medskip\ni
The \nct{dimension} of $G$, denoted $\dim(G)$, is then defined
as the least positive integer $t$ for which $G$ has a realizer of
cardinality $t$.

An intuitive formulation for condition $(*)$ is as follows: For every
vertex $v$ and edge $e$ with $v\not\in e$ the vertex has to get over
the edge in at least one of the orders of a realizer.
%%%%%%%%%%%%%%%%%%%%%%%%%%%%%%%%%%%%%%%%%%%%%%%%%%%%%%%%%%%
%%
\PsFigCap{100}{g-dim}%%
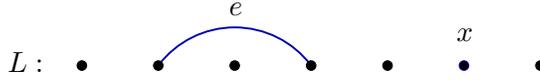
{Vertex $x$ is over edge $e$ in $L$.}
%%
%%%%%%%%%%%%%%%%%%%%%%%%%%%%%%%%%%%%%%%%%%%%%%%%%%%%%%%%%%
All the above results about  dimension of incidence orders of
graphs carry over to this notion of dimension. 
Actually, the two concepts are almost identical:
\Bitem The dimension $\dim(G)$ of a graph equals the interval
dimension of its incidence order $P_G$. In particular
$\dim(G)\leq\dim(P_G)\leq\dim(G)+1$ and $\dim(G)=\dim(P_G)$
if $G$ has no vertices of degree 1 (see \cite{ft-ppg-00}).
\medskip

\ni
Let $L$ and $L'$ be linear orders on a finite set $X$.  We
say that $L'$ is the \term{reverse} of $L$ and write $L'=\rev{L}$ if
$x<y$ in $L$ if and only if $x>y$ in $L$ for all $x,y\in X$.

\begin{defin}
For an integer $t\ge2$, we say that the dimension
of a graph is at most $\below{t}$ if it
has a realizer of the form $\{L_1,L_2,\dots,L_t\}$ with $L_{t}=\rev{L_{t-1}}$.
Similarly, the dimension  is at most $\bbelow{t}$ if it
has a realizer of the form $\{L_1,L_2,\dots,L_t\}$ with $L_{t}=\rev{L_{t-2}}$
and $L_{t-1}=\rev{L_{t-3}}$.
\end{defin}

One of the motivations for introducing this refined version of dimension
was the following theorem proven in~\cite{ft-ppg-00}.
Again, Schnyder woods are the main
ingredient to its proof.

\begin{theorem}\label{thm:outerp}
A graph $G$ is outerplanar iff it has dimension at most $\between{2}{3}$.
\end{theorem}

There are some results concerning the extremal problem of maximizing the
number of edges of a graph of bounded dimension. The first results
from~\cite{ft-ppg-00} only where asymptotic. Felsner~\cite{f-ergd-06}
has obtained sharp bounds:

\Bitem A graph of dimension  $\between{3}{4}$ has at most
$\lfloor\frac{1}{4}n^2 + n -2\rfloor$ edges.

\Bitem A graph of dimension  $\bbetween{3}{4}$ has at most
$\frac{1}{4}n^2 + 5n$ edges.
\medskip

\ni
These bounds easily translate into bounds for the number of characteristic
points of rank~1 on orthogonal surfaces in $\RR^4$ which are generated
by an antichain $V$ with the additional property that certain
pairs of coordinate-orders are reverse to each other.

% *************************************************************
\section{Higher dimensional orthogonal surfaces}\label{sec:hdim}
%*************************************************************

%------------------------------ NON-GENERIC -------------------------------

%%\subsection{Non-generic Surfaces in Higher Dimensions}
%%\label{both}

\subsection{Degeneracies}

On a three-dimensional orthogonal surface, there are three types of
characteristic points: local minima, saddle points, and local
maxima. This classification implies a geometric
\term{rank-function} on the set of characteristic points.

In general, we aim for a combinatorial counterpart for this
concept. In dimension three, points of different
geometric rank are generated in different ways. In the generic case,
the geometric rank of a point coincides with the number of minima
below it. 

\begin{defin}\label{def:MGS} Let $g \in S_V$ be a generated point. A
\term{generating set} for $g$ is a set $G \subset D_g$ such that
$\bigvee(G) = g$. A generating set $G$ is \term{minimal} if $\bigvee(G
\setminus \{v\}) < g$ for all $v \in G$.
\end{defin}

In the generic case, every characteristic point $p$ has a unique
(minimal) generating set, namely $D_p$, so we can simply define the
rank as $r(p) = |D_p| -1$. However, in general, there can be several
minimal generating sets, and they can have different
cardinalities, as illustrated in Figure \ref{fig:degenerate3}.
%
%%%%%%%%%%%%%%%%%%%%%%%%%%%%%%%%%%%%%%%%%%%%%%%%%%%%%%%%%%%
%%
\PsFigCap{40}{degenerate3}%%
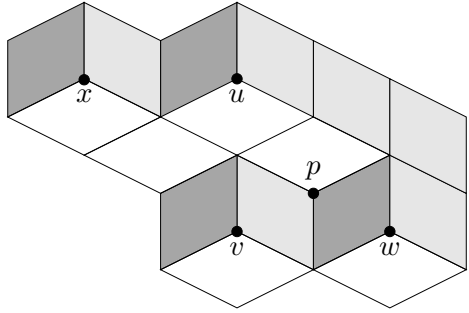
{Point $p$ has minimal generating sets $\{u,v,w\}$ and $\{x,w\}$.}
%%
%%%%%%%%%%%%%%%%%%%%%%%%%%%%%%%%%%%%%%%%%%%%%%%%%%%%%%%%%%
%
The following Lemma shows that such an undesirable situation can be
recognized by a specific pattern in the sets of tight coordinates.

\begin{lemma}\label{lem:pattern} If there is a generated point $g \in
S_V$ with two minimal generating sets of different size, then there
are three minima $u, v, w \in D_g$ and two coordinates $i$ and $j$
such that if we restrict the characteristic vectors $t_g(.)$ of
$T_g(.)$ to positions $i$ and $j$, then we have the following pattern:
\begin{center}\rm
\begin{tabular}{lll} &$i$ &$j$ \\ $t_g(x)$ &1 &1\\ $t_g(u)$ &0 &1\\
$t_g(v)$ &1 &0\\
\end{tabular}
\end{center} In other words, $u_i < v_i = x_i = g_i$ and $v_j < u_j =
x_j = g_j$
\end{lemma}

\Proof For every minimal generating set $G$ and every $v \in
G$, there is some coordinate $i$ such that $w_i < v_i$ for all $w \in
G \setminus \{v\}$, i.e., $v$ is the only minimum
contributing $i$. We call $i$ the \emph{private coordinate} of $v$ in
$G$.

Let $G, G'$ be minimal generating sets for $g$ such that $|G| >
|G'|$. By the pigeon-hole principle, there is some $x \in G'$
covering private coordinates from at least two minima $u,v \in G$.
\hbox to 1em{\hfil}\qed

The converse is almost true. If $g$ is a generated point with the
pattern from Lemma~\ref{lem:pattern} and if $V$ is suspended, then there is a
characteristic point $p$ with minimal generating sets of different
sizes. Such a point $p$ can be reached from $g$ by increasing all
coordinates except the two involved in the pattern. Since the surface is
suspended, each increase is bounded. The point finally reached is
contained in a flat of each color, hence, it is a characteristic point.

\begin{defin}
The antichain $V$ and the corresponding surface $S_V$ are called
\term{degenerate} if there is a characteristic point with the pattern
shown in Lemma \ref{lem:pattern}.
Otherwise, $V$ and $S_V$ are \term{non-degenerate}.
\end{defin}

If $V$ is non-degenerate, then all minimal generating sets of a
characteristic point have the same cardinality. In this case, we
define the \term{rank} of a characteristic point as the size of a minimal
generating set minus one. Minima have rank 0 and  maxima rank $d-1$.

In the 3-dimensional case this definition does not only classify surfaces
with a degenerate vertex as in Figure \ref{fig:degenerate} as `bad', but also
some surfaces which support a proper planar graph, as in Figure
\ref{fig:degenerate3}.

A strong degeneracy is when two different $i$-flats intersect in their
boundaries, as in Figure \ref{fig:degenerate}. From the 3-dimensional
examples it seems plausible that degeneracies which are not strong
could be removed by perturbing flats until different $i$-flats have different
$i$-values, while the cp-order remains the same. The following
example shows that this is not always possible.

We consider a weakly degenerate surface generated by four
minima $a, b, c, d$.  Different $i$-flats have different
$i$-values, hence the characteristic point $p = (2,2,2,2)$ is contained
in exactly four flats. However, since pairs of minima below $p$ share
coordinates (and lie on common flats) in a cyclic structure,
it is not possible to perturb these flats and remove the degeneracy.

$$
a=(2, 2, 1, 1),\quad
b=(1, 1, 2, 2),\quad
c=(2, 1, 2, 1),\quad
d=(1, 2, 1, 2)\quad
$$

We provide the details proving that $a$ and $c$ are contained in a common
1-flat. By the definition of a flat, it is sufficient to find a point
$q \in S_V$ such that $q \rhd_1 a$ and $q \rhd_1 c$. The point $q =
(2,2+\epsilon, 2+\epsilon, 1+\epsilon)$ has the required
properties.

%%%%%%%%%%%%%%%%%%%%%%%%%%%%%%%%%%%%%%%%%%%%%%%%%%%%%%%%%%%%%%%%%%%
%%% Generated vs. characteristic points
%%%%%%%%%%%%%%%%%%%%%%%%%%%%%%%%%%%%%%%%%%%%%%%%%%%%%%%%%%%%%%%%%%%

\subsection{Generated versus characteristic points}

In the following, we assume that $V$ is non-degenerate. 
In this case, we have a combinatorial criterion for characteristic
points:

\begin{prop}\label{prop:char} A generated point $p$ is characteristic
if and only if there are no minima $u, v \in D_p$ such that $T_p(u)
\subset T_p(v)$.
\end{prop}

We will prove this proposition in five steps. Lemmas
\ref{lem:flatMemb1} and \ref{cor:flatMemb2} yield a combinatorial
criterion for the containment of a point $p$ in a given flat $F$.

Lemmas \ref{lem:subset1} and \ref{lem:subset2} establish the
connection between the subset-criterion and 
flat-containment. Finally, Lemma \ref{lem:degOrNonChar} shows that if
$V$ is non-degenerate and $p \in S_V$ is a characteristic point
with $v \in D_p$, then $p_i = v_i \Longleftrightarrow p \in F_i(v)$.

%%%%%%%%%%%%%%%%%%%%%%%%%%%%%%%%%%%%%%%%%%%%%%%%%%%%%%%%%%%%%%%%%%
%%%%%% Lemma: Flat-membership
%%%%%%%%%%%%%%%%%%%%%%%%%%%%%%%%%%%%%%%%%%%%%%%%%%%%%%%%%%%%%%%%%%

\begin{lemma}\label{lem:flatMemb1} 
Let $p \in S_V$. Let $F$ be some
$i$-flat of $S_V$. Then $p \in F$ if and only if there is a $v \in
V_{F}$ and a $q \in S_V$ such that $q \rhd_i v$ and $v \leq p \leq q$.
\end{lemma}

\Proof ``$\Leftarrow$": 
If $p \rhd_i v$, then by definition $p
\in U_i(v) \subset F_i(v)$. Otherwise, $p \neq q$ and there is a $q' \in [p,q]
\subset [v,q] \subset S_V$, such that $q'$ is arbitrarily close to $p$
and $q' \rhd_i v$. Hence $p$ is in the closure of
$U_i(v)$ and thus, $p \in F_i(v)$.

``$\Rightarrow$": If $p$ is in the upper part $F^u$ of $F$, there is
nothing to show , because then $p \rhd_i v$ for some $v \in V_{F}$
and $q = p$ is a good choice for $q$.

Now assume that $p \in F \setminus F^u$. Since $p$ is in the
closure of $F^u$, for every $\epsilon > 0$, there is a $q \in F^u$
such that $|p - q| < \epsilon$. In particular, $|p_j - q_j| <
\epsilon$ for all coordinates $j$.

\Claim 1. There is a $v \in V_{F}$ such that $p \geq v$.
\smallskip

\ni
Let $q \in F^u$ be close to $p$. By definition of $F^u$, there
must be a $v \in V_{F}$ such that $q \rhd_i v$. Suppose $p \ngeq
v$. Then there is a coordinate $j \neq i$ such that $p_j < v_j <
q_j$. This is a contradiction to $|p_j - q_j| < \epsilon $ for
$\epsilon < |p_j - v_j|$.  In short, if $q$ and $p$ are close enough,
then $q\rhd_i v$ implies $p \geq v$.

\Claim 2. There is $q \in F^u$ with $q \geq p$.
\smallskip

\ni
We go for a contradiction and assume 
that there is none. It follows that all $x \in \mathbb{R}^d$ with $x \rhd_i v$
and $x \geq p$ are not on $S_V$. Hence, for every such $x$ there is an
\term{obstructor} $w\in V$ with  $w \lhd x$.

Let ${\bf\gamma}$ be the vector with ${\bf\gamma}_j=1$ iff $v_j=p_j$
and ${\bf\gamma}_j=0$ iff $v_j < p_j$ and ${\bf\gamma}^i$
be obtained from ${\bf\gamma}$ by changing the value of coordinate $i$ to $0$.
Consider the sequence 
$x_n = p + \frac{1}{n}{\bf\gamma}^i$ converging to $p$. If each $x_n$
has an obstructor then there has to be a simultaneous obstructor
$w$ for all elements of the sequence. From $w \lhd x_n$ for all $n$
we obtain that if $j$ is a coordinate with  ${\bf\gamma}^i_j=1$, then
$w_j \leq p_j=v_j$, and if  ${\bf\gamma}^i_j=0$, then $w_j < p_j$.
Hence, for all $j$ either $w_j \leq
v_j$ or $w_j < p_j$.

Let $q \in U_i(v)$ be close to $p$. Since $q$ is not obstructed by $w$
there is a coordinate $j$ with $q_j \leq w_j$. 
Since $q \rhd_i v$ and because $w_j > v_j$ implies $p_j > w_j$
(previous paragraph), we have 
$
v_j < q_j \leq w_j < p_j.
$
This is is impossible if $|p_j - q_j| < \epsilon $ and
$\epsilon < |p_j - w_j|$. It follows that all points in $F^u$ close enough to $p$ are
obstructed by $w$ and, hence, $p\not\in F$. This contradiction
completes the proof.
\qed

Let $p \in S_V$. If $v \in D_p$ is a minimum proving that $p \in F_i(v)$ in
the sense of Lemma \ref{lem:flatMemb1}, i.e. there is a $q \in S_V$
such that $v \leq p \leq q$ and $q \rhd_i v$, then we call $v$ an
\term{$i$-witness} for $p$. 

Obviously, if $v$ is an $i$-witness for $p$, then $p \in F_i(v)$. The
reverse is in general not true.

%%%%%%%%%%%%%%

\begin{corollary}\label{cor:flatMemb2} 
Given $p \in S_V$ and  $v \in D_p$ with
$p_i = v_i$, then $v$ is an $i$-witness for $p$ if and only if there
is no minimum $w \in V$ such that $w_i < v_i = p_i$ and for all $j
\neq i$ either $w_j \leq v_j$ or $w_j < p_j$.
\end{corollary}

\begin{proof} 
This follows from the proof of Claim 2 in the previous
lemma.
\end{proof}

%%%%%%%%%%%%%%

\begin{lemma} \label{lem:subset1} Let $p \in S_V$, $u, v \in D_p$,
$T_p(u) \subset T_p(v)$, $i \in T_p(v) \setminus T_p(u)$. Then $v$ is
no $i$-witness for $p$.
\end{lemma}
\begin{proof} Assume otherwise, and let $q \rhd_i v$ and $q \geq p$. We show that this implies $q \rhd
u$. This is a contradiction to $q \in S_V$. We have to check $q_j >
u_j$ for all $j \in \{1, \dots, d\}$:

For all $j \notin T_p(u)$, we have $q_j \geq p_j > u_j$. This includes
$j = i$. For all $j \in T_p(u)$, $j \neq i$, we have $q_j > v_j$,
because $q \rhd_i v$, and $v_j = u_j = p_j$, because $T_p(u) \subset
T_p(v)$.
\end{proof}

\begin{lemma}\label{lem:subset2} Let $p \in S_V$, $v \in D_p$, $v_i =
p_i$, and assume $v$ is not an $i$-witness for $p$. Then there is a minimum $u \in D_p$ such that $T_p(u)
\subset T_p(v)$ and  $i \in T_p(v) \setminus T_p(u)$.
\end{lemma}

\begin{proof} By Corollary~\ref{cor:flatMemb2}, there is a $u \in V$ such
that $u_i < v_i = p_i$ and for all $j \neq i$: $ u_j \leq v_j$ or $u_j
< p_j$. This implies that there is no coordinate $k$
such that $v_k < u_k = p_k$. Therefore, $T_p(u) \subseteq T_p(v)$ and
since $i \notin T_p(u)$ even  $T_p(u) \subset T_p(v)$ .
\end{proof}

\begin{lemma}\label{lem:degOrNonChar} Let $p \in S_V$, $v \in D_p$,
$v_i = p_i$, but assume $v$ is not  an $i$-witness for $p$.  Then either $p$ is not a
characteristic point or $V$ is degenerate.
\end{lemma}

\begin{proof} Since $v$ is not an $i$-witness for $p$, there is a $u
  \in D_p$ such that $u_i < v_i = p_i$ and $T_p(u) \subset T_p(v)$ by Lemma \ref{lem:subset2}.

Assume $p$ is characteristic. Then $p$ is contained in some $i$-flat,
hence there must be some $i$-witness $w$
for $p$, and $T_p(u) \nsubseteq T_p(w)$ by Lemma
\ref{lem:subset1}. Therefore, there is a coordinate $j \neq i$, $j \in
T_p(u) \setminus T_p(w)$. This results in the following pattern.
\begin{center}
\begin{tabular}{llllll} & & &$i$ &$j$ \\ $t_p(v)$ &= &$(\dots$ &1 &1
&$\dots)$\\ $t_p(u)$ &= &$(\dots$ &0 &1 &$\dots)$\\ $t_p(w)$ &=
&$(\dots$ &1 &0 &$\dots)$\\
\end{tabular}
\end{center}
This shows that $V$ is degenerate.
\end{proof}

Now we can complete the proof of Proposition \ref{prop:char}:

\begin{proof} If $p$ is characteristic and $v \in D_p$,  then by Lemma
\ref{lem:degOrNonChar}, $v$ is an $i$-witness for $p$ if and only if $p_i = v_i$. By
Lemma \ref{lem:subset1}, there can be no minima $u,v \in D_p$ such
that $T_p(u) \subset T_p(v)$.

Conversely let $p$ be generated, every coordinate of $p$ is
covered by some minimum. If there are no $u,v \in D_p$ such that
$T_p(u) \subset T_p(v)$, then Lemma \ref{lem:subset2} implies that $p$
is contained in every flat-type, hence, $p$ is characteristic.
\end{proof}

%%%%%%%%%%%%%%%%%%%%%%%%%%%%%%%%%%%%%%%%%%%%%%%%%%%%%
%%%%%%%% End proof of combinatorial definition
%%%%%%%%%%%%%%%%%%%%%%%%%%%%%%%%%%%%%%%%%%%%%%%%%%%%%%%%

Proposition \ref{prop:char} implies that the minimal generating sets
and the down-set $D_p$ of a characteristic point $p$ have a very
special structure:

\begin{corollary} Let $V$ be non-degenerate. A generated point $p \in
S_V$ is characteristic of rank $k$ if and only if there is a partition
$P_1, \dots, P_{k+1}$ of $D_p$ such that $G \subset D_p$ is a minimal
generating set for $p$ if and only if $|G \cap P_i| = 1$ for all $i =
1, \dots, k+1$.
\end{corollary}

Two minima $u, v \in D_p$ belong to the same part $P_i$ if and only if
$T_p(u) = T_p(w)$.

\begin{corollary} Let $V$ be non-degenerate and $p\in S_V$ be a
  characteristic point. Then every minimum $v \in D_p$ is contained in
  some minimal generating set of $p$.
\end{corollary}

Observe that the corollary does not yield a criterion to
distinguish characteristic points from generated points. 
There exist non-characteristic points such that every minimum below is
contained in a minimal generating set, as in Figure
\ref{fig:degNonChar}.

%%%%%%%%%%%%%%%%%%%%%%%%%%%%%%%%%%%%%%%%%%%%%%%%%%%%%%%%%%%

\PsFigCap{30}{degNonChar} {$p$ has minimal generating sets $\{u, t\}$ and
$\{v, s\}$}

%%%%%%%%%%%%%%%%%%%%%%%%%%%%%%%%%%%%%%%%%%%%%%%%%%%%%%%%%%

\subsection{Syzygy-points and characteristic points}

Algebraists use orthogonal surfaces  as a
tool to obtain resolutions for monomial ideals.
They are specially interested in the syzygy-points
of a surface. In this subsection we discuss the relation
between syzygy-points and characteristic points.

For every point $p \in S_V$, simplicial complex $\Delta_p$
on the set $\{1, \dots, d\}$ is define by:

$$
 I \in \Delta_p 
\quad\Longleftrightarrow\quad
 p + \sum \limits_{i \in I} \epsilon e_i \in S_V
\text{\quad for some $\epsilon > 0$.}
$$

A point $p$ is a \emph{syzygy-point} if $\Delta_p$ has non-trivial
homology.

\begin{lemma}\label{lem:syz} If $p$ is a syzygy-point, then $p$ is
characteristic.
\end{lemma}

\begin{proof} Suppose $p$ not characteristic. Then there is some $i \in
\{1, \dots, d\}$ such that $p$ is not contained in any $i$-flat. This
implies that there is some $\epsilon >0$ such that $p + \epsilon e_i
\in S_V$.

The claim is that in this case $I \in \Delta_p$ implies $I \cup \{i\}
\in \Delta_p$.  Let $I \in \Delta_p$ with $i \notin I$, let $p' = p +
\sum_{j\in I} \epsilon\, e_j$. This point $p'$ is on $S_V$ by
definition of $\Delta_p$ and $p'$ is not contained in any $i$-flat if
$\epsilon$ is small enough. This implies $p' + \epsilon\, e_i = p +
\sum_{j \in I \cup \{i\} } \epsilon\, e_j \in S_V$. Therefore $I \cup
\{i\} \in \Delta_p$.

It follows that every maximal simplex of $\Delta_p$ contains $i$, therefore,
$\Delta_p$ is contractible. Hence, $p$ is not a syzygy-point.
\end{proof}

For $d=3$, the converse is also true: every characteristic point is a
syzygy-point. This can be verified by considering the types of points
shown in Figure~\ref{fig:typen}.

For $d \geq 4$, not all characteristic points are syzygy-points. We
have an example of a characteristic point $p$ that is no
syzygy. $\Delta_p$ is simply a path on four vertices, hence
contractible.
\medskip

\ni
{\bf Example: A characteristic point which is not syzygy}
\smallskip

\ni
The point $p$ is generated by three minima with the following
coordinates:

\begin{center}
\begin{tabular}{lll} $u$ &= &$(2,2,1,1)$\\ $v$ &= &$(2,1,2,1)$\\ $w$
&= &$(1,2,1,2)$\\
\end{tabular}
\end{center}

Let $p = u \vee v \vee w = (2,2,2,2)$. Clearly, $p \in S_V$, because
$p$ shares some coordinate with every minimum below.

We first check that $p$ is indeed a characteristic point, using Lemma
\ref{lem:flatMemb1}. For every $i$, we provide a minimum $m \leq p$ and
a point $q \rhd_i m$ such that $m \leq p \leq q$, thereby proving $p
\in F_i(m)$.

Each of the $q's$ in the following table shares one coordinate with
each of the three minima $u, v, w$. This ensures that $q \in S_V$.

\begin{center}
\begin{tabular}{lll} $p \in F_1(u)$: &$q = (2, 3, 2, 2) \in S_V$, &$q
\rhd_1 u$ and $u \leq p \leq q$\\ $p \in F_2(u)$: &$q = (3, 2, 2, 2)
\in S_V$, &$q \rhd_2 u$ and $u \leq p \leq q$\\ $p \in F_3(v)$: &$q =
(3, 2, 2, 2) \in S_V$, &$q \rhd_3 v$ and $v \leq p \leq q$\\ $p \in
F_4(w)$: &$q = (2, 3, 2, 2) \in S_V$, &$q \rhd_4 w$ and $w \leq p \leq
q$\\
\end{tabular}
\end{center}

Now we examine the simplicial complex $\Delta_p$:

Let $\epsilon > 0$. For all $i = 1, 2, 3, 4$, we have $p + \epsilon
e_i \in S_V$, because $p$ shares two coordinates with every
minimum. This implies that $\Delta_p$ has the vertices 1, 2, 3, 4.

\begin{center}
\begin{tabular}{ll} $p + \epsilon e_1 + \epsilon e_2 \rhd u
\Longrightarrow \{1, 2\} \notin \Delta_p$ &$p + \epsilon e_1 +
\epsilon e_3 \rhd v \Longrightarrow \{1, 3\} \notin \Delta_p$ \\ $p +
\epsilon e_1 + \epsilon e_4 \in S_V \Longrightarrow \{1, 4\} \in
\Delta_p $ &$p + \epsilon e_2 + \epsilon e_3 \in S_V \Longrightarrow
\{2, 3\} \in \Delta_p $ \\ $p + \epsilon e_2 + \epsilon e_4 \rhd w
\Longrightarrow \{2, 4\} \notin \Delta_p$ &$p + \epsilon e_3 +
\epsilon e_4 \in S_V \Longrightarrow \{3, 4\} \in \Delta_p $ \\
\end{tabular}
\end{center}

The complex $\Delta_p$ is a simple path on the vertices $1, 2, 3,
4$, so $\Delta_p$ is contractible and hence has trivial homology. This
proves that $p$ is not a syzygy point.

%%%%%%%%%%%%%%%%%%%%%%%%%%%%%%%%%%%%%%%%%%%%%%%%%%%%%%%%%
%%%%%%%% Rigidity
%%%%%%%%%%%%%%%%%%%%%%%%%%%%%%%%%%%%%%%%%%%%%%%%%%%%%%%%%

\subsection{Rigidity}

Recall from the 3-dimensional case, that rigidity forces
that a characteristic point of rank 1 can only dominate exactly two
minima. An alternative
formulation is that a characteristic point of rank 1 must not dominate
another point of the same rank. This second condition can be
generalized for characteristic points of arbitrary rank:

\begin{defin}
  An orthogonal surface $S_V$ is called \term{rigid} if and only if
  the characteristic points of every rank are an antichain in the cp-order.
\end{defin}

In other words, $V$ is rigid if and only if the cp-order of
$V$ is \term{graded}. For face lattices of
polytopes, this is a necessary condition. In particular, two faces of
a polytope are comparable if and only if one is contained in the
other, and this implies that they have different dimension.

For the three-dimensional case, rigidity is sufficient to ensure that
the dominance order on characteristic points is indeed isomorphic to
the face lattice of a 3-polytope (minus one facet), as we discussed in
Section \ref{sec:3dim}.

However, in dimension four, this is no longer true. There are examples
showing that there remain rather substantial differences between
cp-orders of rigid orthogonal surfaces 
and face lattices of polytopes in general:

\begin{itemize}

\item In a (face) lattice, any two elements have a unique join. There
  are rigid cp-orders that violate this condition and, hence, are no
  lattices (see Subsection \ref{subsec:noLattice}).

\item Even if the cp-order is a
  lattice it may have intervals of height 2 which are no
  quadrilaterals (see Subsection \ref{subsec:noDiamond}).
  This is impossible for face lattices of polytopes (face lattices of polytopes
  have the \emph{diamond-property}).
\end{itemize}

\ni

\begin{problem}
Identify further properties of cp-orders of (rigid) orthgonal surfaces.
\end{problem}

\subsubsection{A rigid flat without the lattice-property}
\label{subsec:noLattice}

Figure~\ref{fig:noLattice} shows one flat $F$ of an orthogonal surface in
dimension 4. The surface is rigid. It is generated by four internal
minima $v,w,s,t$ with coordinates
$$
v=(3, 1, 2, 3),\quad
w=(1, 3, 1, 3),\quad
s=(4, 2, 3, 1),\quad
t=(2, 4, 4, 2)\quad
$$
together with the four suspensions $X,Y,Z,T$. The flat $F$
is the flat $F_4(v) = F_4(w)$ 
In the figure the two minima $v$ and $w$ are marked black,
characteristic points of rank 1 (edges) are marked white, points of
rank 2 (2-faces) are marked blue and maxima are marked green.

Note that the boundary of $F$ consists of a \emph{lower staircase}
containing $v, w$ as well as some characteristic points of rank 1 and 2 and an \emph{upper
staircase} containing the edges $(v,s)$ and $(w,t)$ and all maxima of $F$.
The maximum labeled $M$ is generated by $\{v,w\},s,t,Y$ (this is
to be read as: $M$ is minimally generated by the vertices $s,t,Y$ together
with either $v$ or $w$). 

Consider the interval $[v,M]$ in the dominance order. It is
shown in the right part of the figure. All characteristic points in that interval 
are contained in $F$. Note that $[v,M]$ only contains two edges
$(v,s)$ and $(v,w)$ and two 2-faces generated by $\{v,w\},s,t$ and
$\{v,w\},s,Y$ respectively. Both edges are comparable to both 2-faces in
the dominance order. This shows that the cp-order of the surface 
is not a lattice. 

%%%%%%%%%%%%%%%%%%%%%%%%%%%%%%%%%%%%%%%%%%%%%%%%%%%
%%
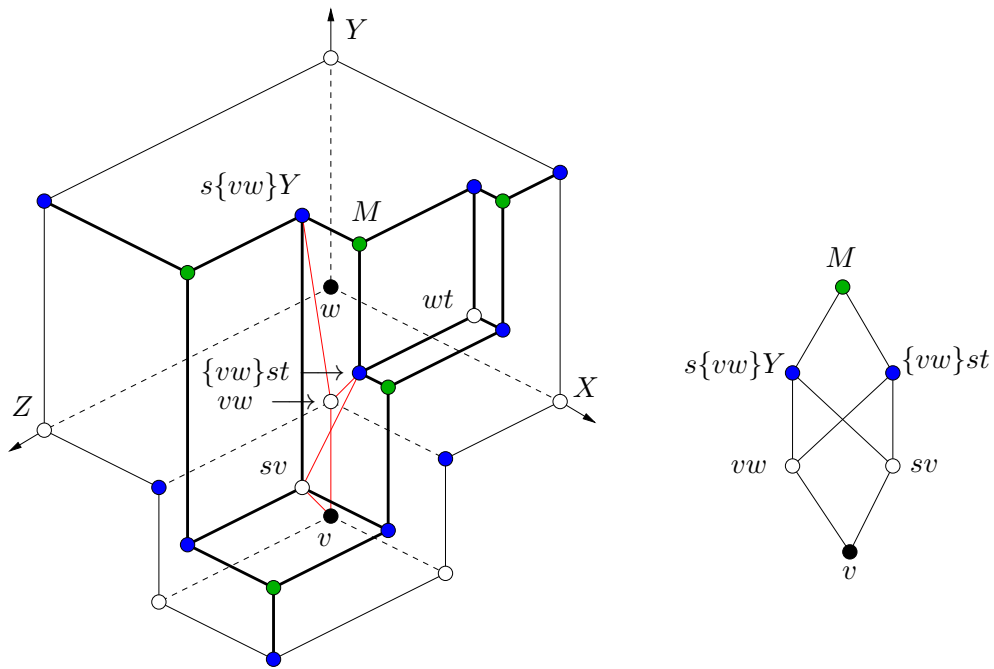
\PsFigCap{60}{noLattice}{A rigid flat violating the Lattice-property}
%
%%
%%%%%%%%%%%%%%%%%%%%%%%%%%%%%%%%%%%%%%%%%%%%%%%%%%%

\subsubsection{A rigid flat without the diamond property}
\label{subsec:noDiamond}

Figure~\ref{fig:noDiamond} shows a flat $F$ of an orthogonal surface
in dimension 4. The surface is rigid. It is generated by the four
suspensions $X,Y,Z,T$ together with six internal minima $x,u,v,w,s,t$:

$$
x = (3, 3, 3, 3),~
u = (1, 4, 4, 3),~ 
v = (4, 1, 4, 3),~
w = (4, 4, 1, 3),~
s = (2, 5, 5, 1),~
t = (5, 5, 2, 1)
$$

(For better visibility we have used a different set of
coordinates in the figure. The combinatorial structure
of the flat is not affected by this change.)
The flat $F$ contains four minima $x,u,v,w$, they share the last
coordinate, so $F$ is a 4-flat.  Consider the characteristic point
$p=(5,5,5,3)$ of rank 2. It is generated by the minima $\{x,u,v,w\},
s,t$ The interval $[v,p]$ contains three characteristic points of rank 1:
$x\join u$ and $x\join v$ and $x\join w$. They are all located on the lower
staircase.

%%%%%%%%%%%%%%%%%%%%%%%%%%%%%%%%%%%%%%%%%%%%%%%%%%%
%%
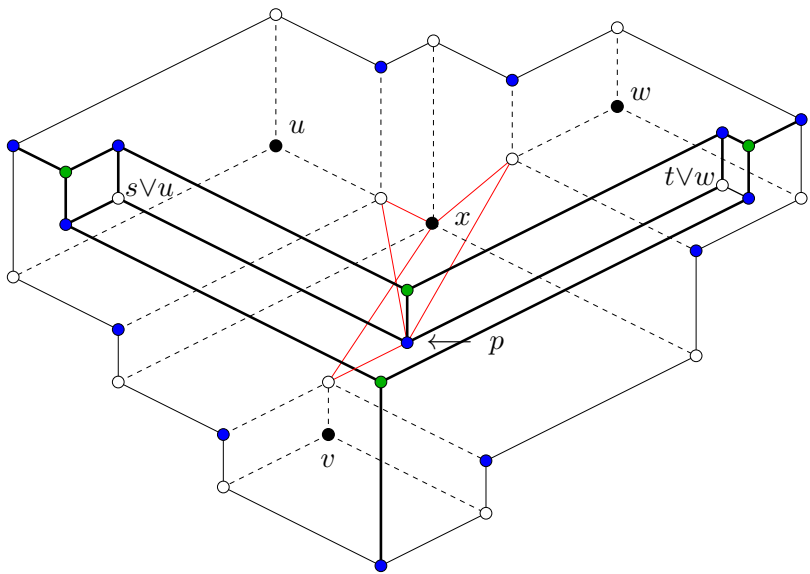
\PsFigCap{55}{noDiamond}{A rigid flat violating the Diamond-property}
%%
%%%%%%%%%%%%%%%%%%%%%%%%%%%%%%%%%%%%%%%%%%%%%%%%%%%

% *************************************************************
\section{Realizability of polytopes}
\label{sec:real}
% *************************************************************

In the previous section we have investigated cp-orders. It became
clear that even non-degenerate rigid surfaces can have cp-orders which
are far from face lattices. In this section we turn the focus to
polytopes and ask whether the face lattice of a given polytope $P$
can be realized on an orthogonal surface.

In Section~\ref{ssec:rot} we have seen realizability criteria
which came from dimension theory of orders. The following
subsection shows a criterion of different guise.
After that we present some families of realizable polytopes.

% *************************************************************
\subsection{Generic surfaces and realizability}
\label{higherDim}

The fact that every Scarf-complex is polytopal immediately raises the
question whether every simplicial polytope is a
Scarf-complex, or \term{realizable}.
However, it is not difficult to find non-realizable  simplicial polytopes.

Every 3-polytope is a Scarf-complex, but already in dimension $d
= 4$, there are large classes of non-realizable simplicial
polytopes. One example, as mentioned in Section \ref{sec:odim}, are
4-polytopes with a skeleton-graph containing the complete graph $K_n$
for $n \geq 13$. There are also smaller examples we will present later
in this section.

A particularly well behaving class of polytopes are \emph{stacked
  polytopes}. These are simplicial polytopes with the minimal number
of faces. Every stacked polytope is realizable, we give a proof for
this later in this section.

For arbitrary simplicial polytopes, there is no complete
characterization  of realizability. However, based on the
combinatorial properties of cp-orders, we have a
necessary criterion. It concerns the number of
incidences between a $k-1$-face and a facet. Validating the criterion
for a specific example only requires counting these incidences.

\begin{prop}[A Realization-Criterion]
  Let $V \subset \mathbb{R}^d$ be a generic suspended antichain, and $p
  \in S_V$ be a characteristic point with $D_p = \{v_1, \dots, v_k\}$.
  For every choice of coordinates $$i_1 \in T_p(v_1), \dots, i_k \in
  T_p(v_k)$$ such that $p_{i_j}>0$ for all $j = 1, \dots, k$, there is
  a maximum $M \in S_V$ such that $M_{i_1} = p_{i_1}, \dots, M_{i_k} =
  p_{i_k}$.
\end{prop}

If $S_V$ is suspended and $p$ is an inner point, i.e. $p_i > 0$ for
all $i$, then this implies that there are at least
$|T_p(v_1)|\cdot|T_p(v_2)|\cdot \ldots\cdot |T_p(v_k)|$ maxima above
$p$. In particular, the cp-lattice of $S_V$ corresponds to a simplicial
polytope where the face with vertices $\{v_1, \dots, v_k\}$ is
contained in at least $\prod_i |T_p(v_i)|$ facets.

\begin{proof} The idea is to start at $p$ and successively augment
  every coordinate $j \notin \{i_1, \dots, i_k\}$ until we reach a
  maximum.  The condition that no $i_j$ is minimal ensures that we do
  not walk into one of the $d$ unbounded flats.

  We walk along the ray $p + \lambda e_j$ for $\lambda >0$. We stop at
  the first point where $p + \lambda e_j \geq u$ for some minimum $u
  \notin D_p$. The point $u$ is unique, because $V$ is generic - there
  is no other minimum with this $j$-coordinate. Furthermore, $u_j \neq
  0$ because $p_j < u_j$. Therfore, we can iterate with the point $p'
  = p \vee u$ with $D_{p'} = D_p \cup \{u\}$ and increase some other
  direction $j' \notin \{i_1, \dots, i_k, j\}$. We can repeat the
  augmentation $d-k$ times. Finally we reach a point that dominates
  $d$ minima and is minimal in no coordinate. This point shares
  exactly one coordinate with every minimum below, so any further step
  in a positive direction would leave the surface. This characterizes
  a maximum $M$. The required properties of $M$ are obvious: $M_{j} =
  p_{j}$ iff $j$ is one of the selected coordinates $i_1,\dots,i_k$.
\end{proof}

We now take a closer look at a special case of this proposition
and its applications. Let $d = 4$ and consider characteristic
points generated by pairs of minima, i.e., edges.

In dimension 4, there are two possible forms for the join of two
minima $u,v$. Either, $u \vee v$ inherits three coordinates from one
of its generators and only one coordinate from the other, so w.l.o.g.
$u \vee v = (u_1, v_2, v_3, v_4)$. The other case is that each
generator contributes two coordinates, so w.l.o.g. $u \vee v = (u_1,
u_2, v_3, v_4)$.

In the first case, $(u,v)$ is the end-point of some orthogonal arc of
$v$, and we call the edge $(u, v)$ an \term{orthogonal edge}. In the
second case, we call it a \term{symmetric edge}. For a symmetric edge
$(u_1, u_2, v_3, v_4)$, the criterion states that there are maxima
$M_{1,3}, M_{1,4}, M_{2,3}, M_{2,4}$ such that $M_{i,j}$ inherits
coordinate $i$ from $u$ and coordinate $j$ from $v$.
In particular, every symmetric edge has to be contained in at  least four facets.

A further important observation is that every inner vertex has exactly
four outgoing orthogonal edges. This implies that in total, there are
exactly $4n -10$ orthogonal edges. All the other edges must be
symmetric. In some cases, we can identify inner edges as symmetric,
because we know that all edges connecting an inner point to a
suspension point are orthogonal edges of the inner point.

These observations yield two useful methods to identify symmetric
edges and thus prove the non-realizability of a complex:

\begin{itemize}
\item Any edge incident to a suspension vertex is orthogonal. Therefore, given two inner vertices $v, w$ such that both are adjacent to \emph{all four suspensions}, we know that the edge $(v,w)$ is symmetric. (``Suspension-criterion")
\item At least  $|E| - (4n-10)$ of the inner edges are symmetric. (``Counting-criterion")
\end{itemize} 

We will consider two different but closely related aspects of the realization problem: 

\begin{itemize}
\item[(A)] Given a simplicial polytope, is there an cp-lattice
  realizing it? \\
  This asks for the realization of a polytopal \emph{sphere}?

\item[(B)] Given a simplicial polytope $P$ with a designated facet $F$, is
  there an orthogonal surface realizing $P$ such that $F$ is the outer
  facet, i.e. the vertices of $F$ are the suspensions?\\
  This asks for the realization of a polytopal \emph{ball}?
\end{itemize}

A polytope is non-realizable in the sense of (A) if and only if it is non-realizable for every choice of a facet $F$ in the sense of (B).
There are polytopes that are realizable in the sense of (A), but not for all
possible choices of $F$. One example (already discussed in 
\cite{bps-mr-98}) is the cyclic polytope $C_4(7)$ . Here is a list of its facets: 
\begin{center}
\begin{tabular}{lllllll}
$[1, 2, 3, 4]$  &$[1, 2, 3, 7]$   &$[1, 2, 4, 5]$   &$[1, 2, 5, 6]$  
&$[1, 2, 6, 7]$ &$[1, 3, 4, 7]$  &$[1, 4, 5, 7]$\\
$[1, 5, 6, 7]$   &$[2, 3, 4, 5]$   &$[2, 3, 5, 6]$ &$[2, 3, 6, 7]$ 
&$[3, 4, 5, 6]$   &$[3, 4, 6, 7]$   &$[4, 5, 6, 7]$\\  
\end{tabular}
\end{center}

The underlying graph is the complete graph $K_7$. This implies that no
matter which facet we choose as the outer facet, every inner vertex
must be adjacent to all four outer vertices, i.e. the suspensions. 
The suspension criterion implies that all edges between inner
vertices are symmetric.

The polytope $C_4(7)$ has two kinds of facets: If we choose one of the
facets listed in the table below as the outer facet, then the
remaining complex is not realizable. There is always an inner edge
that is contained in only three facets.

\begin{center}
\begin{tabular}{llllllll}
  Facet &$[1, 2, 3, 4]$ &$[1, 2, 3, 7]$ &$[1, 2, 6, 7]$
  &$[1, 5, 6, 7]$  &$[2, 3, 4, 5]$ &$[3, 4, 5, 6]$ &$[4, 5, 6, 7]$\\
  Inner E. &$[5, 7]$ &$[4, 6]$ &$[3, 5]$ &$[2, 4]$ &$[1, 6]$ &$[2, 7]$ &$[1, 3]$\\
\end{tabular} 
\end{center}

If we choose any of the seven other facets, the remaining complex is
realizable. 
\medskip

We have used a computer to generate a list of all orthogonal
triangulations on 7, 8 and 9 vertices. This list has been compared
to a list of all simplicial polytopes on these numbers of vertices%
\footnote{This list was compiled by Frank Lutz who also helped with
the computations}.
The comparison yields the following results:

\begin{itemize}
\item All simplicial polytopes on 7 and 8 vertices are realizable.
\item On 9 vertices, there are 116 non-realizable simplicial polytopes
  in the sense of (A). 
\end{itemize} 

\ni Every non-realizable polytopal sphere provides several
non-realizable balls, because we can choose any facet as the outer
facet. The 116 non-realizable polytopes on 9 vertices lead to 2957
non-realizable balls on 9 vertices.  For these examples, we counted
edge-facet incidences and compared to the two realization criteria.
The results:

\begin{itemize}
\item 2141 of the 2957 non-realizable balls violate the
  ``suspension-criterion" (816 do not).
\item 2023 of the 2957 non-realizable balls violate the
 ``counting-criterion" (934 do not).
\end{itemize}

\ni
Together, the two criteria work for 2344 of the balls.
For the remaining 613 balls, counting the
edge-facet incidences is not sufficient to prove that they are
non-realizable.

For some of the difficult examples, new strategies to find symmetric
edges might be sufficient to enable us to use the edge-facet
criterion. For others, the edge-facet criterion is no help, because
every edge belongs to at least 4 facets. For these cases, new
arguments are needed.

% *************************************************************
\subsection{Realizable polytopes}

In this subsection we present classes of polytopes which can be shown to
be realizable by an orthogonal surface. Recall that this means that 
the face lattice of the polytope is a cp-lattice. 
We have already mentioned classes of realizable simplicial polytopes,
like all 4-polytopes with at most 8 vertices and \term {stacked
  polytopes}. These are polytopes that can be constructed from a
simplex by a series of {\itshape stacking operations,} which means that
a facet $F$ is replaced by a vertex $v$ and $d$ new facets. In other words, a
small pyramid with apex $v$ is erected above $F$.

\begin{prop}
Every stacked polytope is realizable on a generic orthogonal surface. 
\end{prop}

\ni
\begin{proof}
A realization can be constructed inductively in the same way as the
stacking, where a stacking operation corresponds to replacing a
maximum $M$ of the orthogonal surface $S_V$ with a vertex $v$ in the
following way: 

Assume $M$ is generated by the $d$ vertices $w_1, \dots, w_d$, where
$w_i$ contributes the $i$th coordinate to $M$. We insert a minimum $v$
with coordinates $v_i = w_{i_i} - \epsilon$, $i = 1, \dots, d$. Obviously, $v \lhd M$. 

For every $i \in \{1, \dots, d\}$, there is a new maximum generated by
$v$ and all $w_j$'s, $j \neq i$ to which $v$ contributes the $i$th
coordinate. These maxima correspond to the $d$ new facets resulting
from the stacking operation.

There are no characteristic points generated by $v$ and any other vertices
besides the  $w_i$'s. Assume there was such point $p > v$. Then $p$ is
strictly greater than $v$ in some coordinate $i$. If $\epsilon$ is
small enough, this implies $w_i \lhd p$, hence $p \notin S_V$. 

\end{proof}

\begin{prop}
All $d$-polytopes on $d+2$ vertices are realizable.
\end{prop}

\ni
{\sl Proof (sketch).}
$d$-polytopes with $d+2$ vertices are completely classified,
\cite{g-cp-03}.
For every $d$, there are $\lfloor \frac{d}{2} \rfloor$ combinatorial
types of simplicial $d$-polytopes with $d+2$ vertices. 
Non-simplicial $d$-polytopes on $d+2$ vertices are pyramids over some
$d-1$-polytope on $d+1$ vertices.

A proof of  the two following facts can be found in \cite{k-hdos-06}:
\begin{itemize}
\item If a $d-1$-polytope $P$ is a realizable on an orthogonal surface of
dimension $d-1$, then the pyramid $pyr(P)$ is realizable in dimension
$d$. 
\item There are $\lfloor \frac{d}{2} \rfloor$ combinatorially different
orthogonal triangulations in dimension $d$ on $d+2$ vertices.
\end{itemize}

\ni
The proposition follows by induction from these facts together with the
realizability of all 3-polytopes. 
\qed
\bigskip

\ni
Products of a polytope with an edge and more generally products with
paths ({\itshape sequences}) preserve realizability.
The construction is detailed in \cite{k-hdos-06}, here we
only indicate the ideas:

\begin{proposition}
Let $P$ be a $d-1$-polytope. If $P$ is a facet of some realizable
$d$-polytope, then the prism over $P$, i.e. the product of $P$ with an
edge,  is also a  realizable $d$-polytope.
\end{proposition}

\ni
{\sl Proof (sketch).}
Assume that $P$ is realized as maximum $M$ in an orthogonal surface of
dimension $d$. The idea is to duplicate the local structure of $M$ with
slightly perturbed coordinates.

If $v$ is a minimum contributing coordinate $i$ to $M$, its
double is $v' = v + 2 \epsilon e_i - \epsilon {\bf 1}$.

It is easy to check that the new vertices leave the old structure
unchanged, i.e. they cannot obstruct any old faces. However,
the join of the set of all new vertices is obstructed by
$M$. Therefore, the counterpart-facet $M'$ of $M$ has no corresponding
point on the surface. $M'$ is the outer facet of the realization,
i.e., the facet of the polytope which is missing in the cp-lattice of
the surface.
\qed

\begin{corollary}
The $d$-cube is realizable.
\end{corollary}

\smallskip

\ni The following more general construction produces a realization of 
a $P$-sequence, i.e., of the
product of a realizable $d-1$-polytope $P$ with a path: The
realization of the product consists of translated copies of
realizations of $P$:
  
Given a $d-1$-realization of $P$ with minima $v_1, \dots, v_n$, let
$v^i_j$ is the copy of vertex $v_j$ in $P^i$. The 
coordinates of $v^i_j$ are $(v_j - i \epsilon , i)$. For a small enough
$\epsilon>0$, the set $V_i = \{v^i_j : j=1\ldots n\}$ is an antichain.
It is also easy to see that a vertex
$v^i_j$ is adjacent to only two vertices outside $V_i$, namely
$v^{i-1}_j$ and $v^{i+1}_j$.

\medskip

It is an interesting question to identify further classes of
realizable polytopes. Since the cyclic polytope $C(n,d)$  is not
realizable for $n$ sufficiently large we are particularly
curious about the following:

\begin{problem}
Are the dual polytopes of cyclic polytopes always realizable?
\end{problem}

% *************************************************************

{\advance\baselineskip-1pt
\bibliography{os-bib}
\bibliographystyle{/homes/combi/felsner/FU-DATA/TexStuff/my-siam}
}

% *************************************************************
\end{document}